\documentclass[a4paper,11pt]{article}

\usepackage[english]{babel}
\usepackage[utf8]{inputenc}
\usepackage{amsmath}
\usepackage{amssymb}
\usepackage{graphicx}
\usepackage[colorinlistoftodos]{todonotes}
\usepackage[a4paper]{geometry}

\usepackage{enumerate}
\newtheorem{theorem}{Theorem}[section]
\newtheorem{proposition}[theorem]{Proposition}

\newtheorem{definition}[theorem]{Definition}
\newtheorem{remark}[theorem]{Remark}

{}

\begin{document}

\newgeometry{left = 3cm, right = 3cm, top = 3cm, bottom = 3cm}

\title{Diagonal stability of a class of discrete-time positive switched systems with delay}

\author{Alexander Aleksandrov \thanks{Faculty of Applied Mathematics and Control Processes, Saint Petersburg State University, Saint Petersburg, 199034, Russia, and Department of Control of Complex Systems, 
ITMO University, Saint Petersburg, 197101, Russia. email: alex43102006@yandex.ru}, Oliver Mason\thanks{Dept. of Mathematics and Statistics/Hamilton Institute, Maynooth University, Maynooth, Co. Kildare, Ireland and Lero, the Irish Software Research Centre. email: oliver.mason@mu.ie}}

\date{\today}
\maketitle 

\begin{center}
Submitted to IET Control Theory and Applications
\end{center}

\begin{abstract}
A class of discrete-time nonlinear positive time-delay switched systems with sector-type nonlinearities is studied. Sufficient conditions for the existence of common and switched diagonal Lyapunov--Krasovskii functionals for this system class are derived; these are expressed as feasibility conditions for systems of linear algebraic inequalities.  Corresponding spectral conditions for the existence of common L--K functionals are also described.  Furthermore, it is shown that the proposed approaches can be applied to discrete-time models of digital filters and neural networks. Finally, a numerical example is given to illustrate the effectiveness of theoretical results.
\end{abstract}

\section{Introduction}\label{sec1}
Positive systems, for which nonnegative initial conditions give rise to nonnegative trajectories \cite{FarRin, ValCol}, are of significant practical importance due to their role in various applications including: population dynamics; consensus problems; and congestion control in transmission control protocol networks \cite{ShoTCP}.  While the fundamental properties of positive linear time-invariant (LTI) systems are now well studied and understood \cite{FarRin}, there is a practical need to extend this theory to more complex, realistic models.  In particular, issues such as nonlinearities, time-delay and switches in system dynamics \cite{ValCol} give rise to novel problems, many of which are not yet resolved.  The issue of time-delay is of particular importance for applications to networked control systems \cite{ZhaoLiuRees09}.  In the current paper, we are concerned with the fundamental question of stability for a class of nonlinear, switched positive systems subject to time-delay.  The model class considered is in discrete time. 

The stability theory of positive systems has several aspects that distinguish it from the corresponding theory for general systems.  Two of the most significant of these are arguably the use of copositive Lyapunov functions \cite{MasSho07} and diagonal Lyapunov functions \cite{Kaz,Lam,AlexMas} in the stability analysis of positive systems.  It is well known, and follows from classical Perron-Frobenius theory \cite{HJ}, that the stability of a positive LTI system is equivalent to: i) the existence of a linear copositive Lyapunov function; and ii) the existence of a diagonal Lyapunov function.  These simple facts have motivated researchers to investigate the existence of corresponding types of Lyapunov functions and Lyapunov--Krasovskii (L--K) functionals for systems subject to time-delay \cite{Lam,AlexMas}, switching \cite{Short,Past,Huang} and nonlinearities \cite{AiT,Shaker}.  Our work here continues in this vein, focusing on the existence of common diagonal and switched diagonal L--K functionals for a class of switched positive systems with time-delay.  Specifically, we shall first present conditions for the existence of common diagonal L--K functionals for the system class that are less conservative than those given in \cite{AlexMas}.  We shall also describe results on the existence of switched diagonal L--K functionals that give conditions in terms of simple algebraic inequalities in contrast to the LMI conditions given in \cite{Shaker} where undelayed systems subject to a related but narrower class of nonlinearities were considered.      

There are several reasons for considering the question of diagonal stability for this system class.  First of all, diagonal functionals are attractive because of their simple structure in which individual states are decoupled.  In many cases, the existence of a diagonal Lyapunov function implies a more robust form of stability, ensuring that the system is stable even when subject to parameter uncertainty \cite{Kaz}.  From a practical point of view, common and switched diagonal Lyapunov functions can be used to stabilize switched linear systems \cite{ZhenFen12}.  Moreover, such functions have proven to be powerful tools in the analysis of neural networks, asynchronous computation and in the so-called large scale systems approach \cite{Kaz}.  Finally, on a theoretical level, it is of course interesting to understand the degree to which properties of basic LTI systems can be extended to more complex system classes and results characterising the existence of diagonal L--K functionals allow us to understand the, somewhat intricate, relationship between the various types of Lyapunov function available.  In this context, the interesting work of \cite{Past} on the links between linear, max-type and diagonal Lyapunov functions for switched positive linear systems is noteworthy; in fact, the work of this latter paper has inspired some of the results to be presented here.  

\subsection{Contributions}\label{sec:con}
The main contributions of this paper are described below.

\begin{itemize}
\item[(i)] We describe less conservative conditions for the existence of common diagonal L--K functionals for the considered classes of subsystems than those described in \cite{AlexMas}.  Equivalent spectral conditions are also described. An explicit numerical example to highlight the reduction in conservatism is given.
\item[(ii)] We provide constructive, readily verifiable conditions for the existence of common and switched L--K functionals.  Furthermore, while sufficient conditions in terms of LMIs have been given previously \cite{Shaker}, our conditions are formulated in terms of linear algebraic inequalities, and can be solved using linear programming.  Links to the spectral radii of matrices associated with the systems are also described.
\item[(iii)] We extend the result of \cite{Past}, which applied to linear switched positive delay-free systems, to a class of nonlinear switched positive systems with delay. In addition, compared with \cite{Past}, we derive conditions for the existence of not only common diagonal L--K functionals, but also switched diagonal L--K functionals.
\end{itemize}

\section{Notation and background}
Throughout the paper $\mathbb{R}^n$ and $\mathbb{R}^{n \times n}$ denote the vector spaces of $n$-tuples of real numbers and of
$n \times n$ matrices respectively.  The notation $\|\cdot\|$ refers to the Euclidean vector norm. For vectors $v \in \mathbb{R}^n$, $v \geq 0$ ($v \leq 0$) means $v_i \geq 0$ ($v_i \leq 0$) for $1 \leq i\leq n$, and $v \gg 0 $ ($v \ll 0 $) means 
$v_i > 0$ ($v_i < 0$) for $1 \leq i \leq n$.  We use the notation $A^T$ for the transpose of a matrix $A$
and $P \succ 0$ ($P \prec 0$) to denote that the matrix $P$ is positive (negative) definite.  Let $\textrm{diag}\left(c_1,\ldots,c_n\right)$
indicate a diagonal matrix with the elements
$c_1,\ldots,c_n$ along the main diagonal.
The identity matrix is denoted by $I$; the dimension will be clear in context.  We say that a matrix $A \in \mathbb{R}^{n \times n}$ is nonnegative if all of its entries are nonnegative.  The \emph{spectral radius} of $A \in \mathbb{R}^{n \times n}$ is denoted by $\rho(A)$.  A matrix $A$ is Schur if all of its eigenvalues have modulus strictly less than 1 ($\rho(A) < 1$); $A$ is Hurwitz if all of its eigenvalues have negative real parts.  A matrix $A$ is Metzler if its off-diagonal entries are all nonnegative.

The following known facts about nonnegative and Metzler matrices are useful for our later results, see \cite{HJ,Kaz}.
\begin{proposition}
\label{prop:Schur} Let $A \in \mathbb{R}^{n\times n}$ be nonnegative.  Then
\begin{itemize}
\item[(i)] $A$ is Schur if and only if there exists some $v \gg 0$ with $A v \ll v$;
\item[(ii)] if $A\theta \ll \theta$ and $A^T d \ll d$, where $\theta\gg 0$, $d\gg 0$, then defining $P = \textrm{diag}\left(d_1/\theta_1,\ldots,d_n/\theta_n\right)$,
we obtain  $A^TPA - P \prec 0$.  
\end{itemize}
\end{proposition}

\begin{proposition}
\label{prop:Metzl} Let $A \in \mathbb{R}^{n\times n}$ be Metzler.  Then
\begin{itemize}
\item[(i)] $A$ is Hurwitz if and only if there exists some $v \gg 0$ with $A v \ll 0$;
\item[(ii)] if $A\theta \ll 0$ and $A^T d \ll 0$, where $\theta\gg 0$, $d\gg 0$, then defining $P = \textrm{diag}\left(d_1/\theta_1,\ldots,d_n/\theta_n\right)$,
we obtain  $A^TP+PA \prec 0$.  
\end{itemize}
\end{proposition}

Our main results in this paper are expressed in terms of systems of algebraic linear inequalities, which are related to the existence of linear copositive Lyapunov functions or the so-called S-property for matrix sets; for background on these and related topics, see \cite{Knorn, ForVal10, Doan, Gowda}.  We here recall some notation and terminology related to these questions.  

Given a finite set of matrices, $M$, the \emph{row selection set} $\mathcal{R}(M)$ of $M$ consists of all matrices formed from the elements of $M$ in the following manner.  The $i$th row of each matrix in $\mathcal{R}(M)$ is the $i$th row of some matrix in $M$ for each $i$, $1 \leq i \leq n$.  Given a set of Metzler matrices $M$, there exists a $v \gg 0$ with $A v \ll 0$ for all $A \in M$ if and only if every matrix in $\mathcal{R}(M)$ is Hurwitz.   The following result is a simple consequence of this fact and has been proven in a distinct manner using \emph{Collatz-Wielandt} sets in \cite{Doan}.

\begin{theorem}
\label{thmcoplin} Let $M$ be a finite set of nonnegative matrices in $\mathbb{R}^{n \times n}_+$.  There exists a vector $v \gg 0$ such that $Av \ll v$ for all $A \in M$ if and only if every matrix in $\mathcal{R}(M)$ is Schur.
\end{theorem}

\section{Statement of the problem}

Consider the switched system
\begin{equation}\label{eq:sys1}
x(k+1)=A_{\sigma(k)} f(x(k))+B_{\sigma(k)} f(x(k-l))   
\end{equation}
and the corresponding family of subsystems

\begin{equation}\label{eq:sys2}
x(k+1)=A_{s} f(x(k))+B_{s} f(x(k-l)), \quad s=1,\ldots,N.   
\end{equation} 
Here $x(k)\in \mathbb{R}^n$;
$A_1, \ldots, A_N$, $B_{1}, \ldots, B_{N}$ are constant matrices; the nonlinearity  
$f: \mathbb{R}^n \rightarrow \mathbb{R}^n$ is continuous and diagonal, meaning: 
$
f(x) = (f_1(x_1),\ldots, f_n(x_n))^T.
$
Furthermore, we assume that each $f_i$ satisfies the following conditions:
\begin{equation}\label{eq:f1}
x_if_i(x_i)>0 \ \ \ \text{for}\ \ x_i\neq 0, 
\end{equation}

\begin{equation}\label{eq:f2}
|f_i(x_i)|\leq |x_i|.            
\end{equation}
%

Finally, $l$ is a positive integer delay.
The switching signal $\sigma$ maps the nonnegative integers into $\{1, \ldots, N\}$ and selects which constituent subsystem is active at each time $k$.  

Systems of this class are closely related to continuous time systems of Persidskii type \cite{Kaz} and are motivated by numerous applications, such as control systems, digital filters, neural networks and iterative numerical methods; see, for example, \cite{Kaz,Phi,Sontag,Liao,EM}.

In what follows we will assume that the matrices
$A_1, \ldots, A_N$, $B_{1}, \ldots, B_{N}$ are nonnegative. Under this assumption, \eqref{eq:sys1} defines a switched positive time-delay system in discrete time, meaning that if initial conditions $x(-l), \ldots, x(0)$ are nonnegative, then $x(k) \geq 0$ for all $k \geq 0$.  For notational convenience, we write $x^{(k)}$ for the augmented state vector $x^{(k)} = \left
( x^T(k), x^T(k-1), \ldots , x^T(k-l) 
\right
)^T$ in $\mathbb{R}^{(l+1)n}$.

We will derive conditions for the existence 
of a common diagonal L--K functional  of the form
\begin{eqnarray}
\nonumber
V(x^{(k)})&=&x^T(k)Px(k)+f^T(x(k-1))Q_1f(x(k-1)) \\
\label{eq:V1}
&+&\ldots+f^T(x(k-l))
Q_l f(x(k-l)) 
\end{eqnarray}
for the family \eqref{eq:sys2}. Here $P,Q_1,\ldots,Q_l$ are positive definite diagonal matrices.
If the family \eqref{eq:sys2} admits such a functional whose differences, $V(x^{(k+1)}) - V(x^{(k)})$ along trajectories of all constituent subsystems are negative, then we say that the system \eqref{eq:sys1} is \emph{diagonally stable}.

\begin{remark}
If the system \eqref{eq:sys1} is diagonally stable, then its zero solution is asymptotically stable for any nonlinearities $f_1(x_1),\ldots,f_n(x_n)$ and for any switching law.
\end{remark}

The problem of diagonal stability of \eqref{eq:sys1} is equivalent to the feasibility of the associated system of LMIs.
This problem has been well investigated, see, for instance, \cite{Kam,Boid}, and various numerical schemes are available to determine LMI feasibility.  However, it should be noted that for systems with parametric uncertainty and for situations where understanding the relationship between the existence of certain types of Lyapunov function and the dynamical properties of the underlying system, the determination of 
simple analytic criteria is crucial.

In \cite{AlexMas}, sufficient conditions for the diagonal stability of the system \eqref{eq:sys1} were derived; these conditions were formulated in terms of the feasibility of auxiliary systems of linear algebraic inequalities.  Formally, the following theorem was proved in \cite{AlexMas}.

\begin{theorem}\label{thm1}
  Let  there exist  
vectors $d \gg 0$, $\theta \gg 0$ such that
\begin{equation}\label{eq:d1}
(A_{s}+B_{r})^T d \ll d, \qquad s,r=1,\ldots ,N,
\end{equation}

\begin{equation}\label{eq:theta1}
(A_{s}+B_{s})\theta  \ll \theta, \qquad s=1,\ldots ,N. 
\end{equation}
Then the system \eqref{eq:sys1} is diagonally stable.
\end{theorem}

\begin{remark}\label{rem1a}
Using Theorem \ref{thmcoplin}, it is possible to state the conditions given in Theorem \ref{thm1} in terms of the spectral radii of sets of matrices associated with the system.  First of all, denote by $M_1$ the set of all matrices of the form $(A_s + B_r)^T$ for $1 \leq s, r \leq N$; similarly, denote by $M_2$ the set of all matrices $A_s+B_s$ for $1 \leq s \leq N$.  Now letting $\mathcal{R}(M_1)$ and $\mathcal{R}(M_2)$ be the sets of all {row selections} from $M_1$, $M_2$ respectively, the theorem can be reformulated as follows.
\end{remark}

\begin{theorem}\label{thm1a}
Let $\rho_1 = \max\{\rho(A) \mid A \in \mathcal{R}(M_1)\}$ and $\rho_2 = \max\{\rho(A) \mid A \in \mathcal{R}(M_2)\}$.  If $\rho_1 < 1$, $\rho_2 < 1$ then the system \eqref{eq:sys1} is diagonally stable.
\end{theorem}

Our first goal here is to relax the conditions for diagonal stability formulated in Theorem \ref{thm1}.  Furthermore, along with conditions of the existence of a common diagonal L--K functional, conditions of the existence of a switched diagonal L--K functional for system \eqref{eq:sys1} will be obtained.

\section{Construction of a common diagonal L--K functional}

In \cite{Past}, an approach was described that allows us to relax conditions for the existence of a common quadratic diagonal Lyapunov function for a family of linear positive delay-free systems. In the present section, we apply the core idea of this approach to derive new conditions for diagonal stability of the nonlinear switched positive time-delay system~\eqref{eq:sys1}.

\begin{theorem}\label{thm:Swit2}
Let there exist numbers $\mu>0$, $\lambda>0$ and 
vectors $d \gg 0$, $\theta \gg 0$ satisfying the
inequalities
\begin{equation}\label{eq:d2}
(A_{s}+B_{r})^T d \leq \mu d, \qquad s,r=1,\ldots ,N,
\end{equation}

\begin{equation}\label{eq:theta2}
(A_{s}+B_{s}) \theta \leq \lambda \theta, \qquad s=1,\ldots ,N, 
\end{equation}
\begin{equation}\label{eq:mul}
\lambda \mu< 1. 
\end{equation}
Then the system \eqref{eq:sys1} is diagonally stable.
\end{theorem}
\textbf{Proof}  We will show how to construct a common L--K functional for the family \eqref{eq:sys2}, given by formula \eqref{eq:V1}.  Set $P=\textrm{diag}\left(d_1/\theta_1,\ldots,d_n/\theta_n\right)$, $d_i$ and $\theta_i$ are the components of the vectors $d$ and $\theta$ respectively.  We will show how to choose a positive definite diagonal $Q$ and a real number $\varepsilon>0$  such that defining $Q_m=Q+(l-m+1)\varepsilon I$, $m=1,\ldots,l$, yields a diagonal L--K functional of the form \eqref{eq:V1} for the family of systems \eqref{eq:sys2}.

Consider the difference $\Delta V = V(x^{(k+1)}) - V(x^{(k)})$ of the functional \eqref{eq:V1} along trajectories of the $s$-th subsystem from family \eqref{eq:sys2} for some $s$ in $\{1,\ldots,N\}$. We obtain
$$
\Delta V \leq W_{1s}(x^{(k)})+W_{2s}(x^{(k)}).
$$
Here
$$
W_{1s}(x^{(k)})=
f^T(x(k))\left(A_s^T P A_s-P+  Q \right)f(x(k))
$$
$$
+2f^T(x(k)) A_s^T P B_sf(x(k-l))
$$
$$
+
f^T(x(k-l))\left(B_s^T P B_s-  Q \right)f(x(k-l)),
$$
$$
W_{2s}(x^{(k)})=\varepsilon\left(l\|f(x(k))\|^2-
\sum_{m=1}^l\|f(x(k-m))\|^2\right).
$$

The function $W_{1s}(x^{(k)})$ is a quadratic form in $f(x(k))$ and $f(x(k-l))$ defined by the matrix
$$
C_s = \left(\begin{array}{c c} 
				 A_s^T P A_s-P+  Q   &  A_s^T P B_s\\
			     B_s^T P A_s  &  B_s^T P B_s-  Q  			\end{array}\right).
			$$

It is straightforward to verify that
$$
C_s \left(\begin{array}{c} 
				 \theta \\
			     \theta   			\end{array}\right)=
			     \left(\begin{array}{c} 
				 A_s^T P (A_s+B_s)\theta-d+Q\theta \\
			     B_s^T P (A_s+B_s)\theta-Q\theta   			\end{array}\right)
			$$
$$
\leq \left(\begin{array}{c} 
				 \lambda A_s^T d -d+Q\theta \\
			     \lambda B_s^Td-Q\theta   			\end{array}\right).
$$

Define $\eta=\max_{r=1,\ldots,N}\left\{B_r^Td\right\}$, where the maximum is taken componentwise. Next, choose $Q$ to be a diagonal positive definite matrix such that $Q\theta=\lambda \eta+\delta e$. Here $\delta$ is a positive parameter, and $e=(1,\ldots,1)^T$.
Then the inequalities
$$
C_s \left(\begin{array}{c} 
				 \theta \\
			     \theta   			\end{array}\right)\leq 
			     \left(\begin{array}{c} 
				 (\lambda\mu-1)d+\delta e \\
			     -\delta e   			\end{array}\right)\ll 0
$$
are valid for sufficiently small values of $\delta$. 

As $\theta\gg 0$ and $C_s$ is Metzler, it follows from Proposition \ref{prop:Metzl} that $C_s$ is a Hurwitz matrix. Moreover, $C_s$ is a symmetric matrix. Therefore, it is negative definite. Hence, one can choose a value of the parameter $\varepsilon$ and a number $\beta>0$ such that the estimate 
$$
\Delta V \leq -\beta\sum_{j=0}^l  \|f(x(k-j))\|^2
$$
holds.

\begin{remark}\label{rem2a} Using Theorem \ref{thmcoplin}  and the notation defined in Remark \ref{rem1a}, we can reformulate the previous result and provide equivalent spectral conditions for diagonal stability.  The following fact is key.
\end{remark}
 
\begin{proposition}\label{prop:equiv}
Let $\rho_1 = \max\{\rho(A) \mid A \in \mathcal{R}(M_1)\}$ and $\rho_2 = \max\{\rho(A) \mid A \in \mathcal{R}(M_2)\}$.  There exist real numbers $\mu > 0$, $\lambda > 0$ and vectors $d \gg 0$, $\theta \gg 0$ in $\mathbb{R}^n$ satisfying \eqref{eq:d2}--\eqref{eq:mul} if and only if $\rho_1\rho_2 < 1$.
\end{proposition}
\textbf{Proof}
First assume that $\rho_1 \rho_2 < 1$.  It follows that there is some $\epsilon  >0$ such that $(\rho_1 + \epsilon)(\rho_2 + \epsilon) < 1$.  From the definitions of $\rho_1$, $\rho_2$, all matrices in the sets $\frac{1}{\rho_1 + \epsilon}\mathcal{R}(M_1)$, $\frac{1}{\rho_2 + \epsilon}\mathcal{R}(M_2)$ are Schur.  Hence by Theorem \ref{thmcoplin} there exist vectors $d \gg 0$ and $\theta \gg 0$ such that $(A_s + B_r)^Td \ll (\rho_1 + \epsilon) d$, $(A_s + B_s) \theta \ll (\rho_2 + \epsilon) \theta$ for $1 \leq s, r \leq N$, and $(\rho_1 + \epsilon)(\rho_2 + \epsilon) < 1$.

Conversely, assume that there exist real numbers $\mu > 0$, $\lambda > 0$ and vectors $d \gg 0$, $\theta \gg 0$ in $\mathbb{R}^n$ satisfying \eqref{eq:d2}--\eqref{eq:mul}.  It follows that for any $\epsilon > 0$, $(A_s + B_r)^T d \ll (\mu+\epsilon)d$, $(A_s + B_s)\theta \ll (\lambda + \epsilon) \theta$ for $1 \leq s, r \leq N$.  Theorem \ref{thmcoplin} then implies that $\rho_1 < \mu+\epsilon$, $\rho_2 < \lambda + \epsilon$.  As $\epsilon > 0$ is arbitrary, it follows that $\rho_1 \leq \mu$, $\rho_2 \leq \lambda$  and hence that $\rho_1 \rho_2 < 1$.  This completes the proof.

The following reformulation of Theorem \ref{thm:Swit2} is now immediate. 

\begin{theorem}\label{thm2a}
Let $\rho_1 = \max\{\rho(A) \mid A \in \mathcal{R}(M_1)\}$ and $\rho_2 = \max\{\rho(A) \mid A \in \mathcal{R}(M_2)\}$.  If $\rho_1 \rho_2 < 1$ then the system \eqref{eq:sys1} is diagonally stable.
\end{theorem}
 
This alternative formulation makes clear how the conditions of Theorem \ref{thm:Swit2} relax those given in Theorem \ref{thm1}.  Moreover, it gives a spectral (albeit non-constructive) method of verifying the conditions of the theorem.  It also shows how Theorem \ref{thm:Swit2} is a direct extension of Theorem 3 in \cite{Past} to the time-delayed, nonlinear system class considered here. 
 
 \begin{remark}
To verify the existence of 
 numbers $\mu>0$, $\lambda>0$ and 
vectors $d \gg 0$, $\theta \gg 0$ satisfying 
inequalities \eqref{eq:d2}--\eqref{eq:mul}, one can also adapt the procedure proposed in \cite{Past}. 
\end{remark}

We note that Theorem \ref{thm:Swit2} can be extended to systems with multiple delays.

Formally, consider the family of subsystems
\begin{eqnarray}
\nonumber x(k+1)=A_{s}f(x(k))+B_{1s}
f(x(k-1))\\
\label{eq:11}
+\ldots+
B_{ls}f(x(k-l)), \quad s=1,\ldots,N, 
\end{eqnarray}
where $A_{s},B_{1s}\ldots,B_{l s}$ are nonnegative matrices for $s=1,\ldots,N$, and the remaining notation is the same as used for \eqref{eq:sys2}.

\begin{theorem}\label{thm:Swit3}
Assume that there exist numbers $\mu>0$, $\lambda>0$ and 
vectors $d \gg 0$, $\theta \gg 0$ such that
$$
(A_{s}+B_{1 r_1}+\ldots+B_{l r_l})^T d \ll \mu d, \quad s,r_1,\ldots,r_l
=1,\ldots,N, 
$$
$$
(A_{s}+B_{1 s}+\ldots+B_{l s}) \theta \ll \lambda \theta , \quad s=1,\ldots,N, 
$$
and inequality \eqref{eq:mul} holds.
Then there exists a common diagonal L--K functional  of the form \eqref{eq:V1} for the family \eqref{eq:11}.
\end{theorem}

The proof of the theorem is essentially the same as that of Theorem \ref{thm:Swit2} although the notation is a little more involved.

\section{Construction of a switched diagonal L--K functional}

One approach to obtain less conservative stability conditions for switched systems than those obtained with the aid of common Lyapunov functions is based on the construction of switched Lyapunov functions \cite{Daaf,Shaker,Liu,zapp}. In particular, in \cite{Daaf}  and \cite{Shaker} 
switched quadratic Lyapunov functions were used for the stability analysis of discrete linear delay-free switched systems and discrete nonlinear delay-free switched $\Phi$-systems respectively.  We will adapt this approach to construct switched diagonal L--K functionals for the system \eqref{eq:sys1}.

Specifically, we will construct a switched diagonal L--K functional of the form
\begin{eqnarray}
\nonumber \widetilde V(k,x^{(k)})=x^T(k)P^{(\sigma(k))}x(k)\\
\nonumber 
+f^T(x(k-1))Q_1f(x(k-1))\\
\nonumber 
+f^T(x(k-2))Q_2f(x(k-2))\\
\label{eq:12}
+\ldots+f^T(x(k-l))
Q_lf(x(k-l)) 
\end{eqnarray}
for the time-delay system~\eqref{eq:sys1}. Here $P^{(1)},\ldots,P^{(N)},Q_1,\ldots,Q_l$ are positive definite diagonal matrices.

In contrast to \cite{Daaf,Shaker}, we will give conditions for the existence of such functionals that are based on systems of linear
algebraic inequalities.

\begin{theorem}\label{thm:Swit4}
Assume that there exist numbers $\mu>0$, $\lambda>0$ and 
vectors $d^{(1)} \gg 0,\ldots,d^{(N)}\gg 0, \theta \gg 0$ such that
\begin{equation}\label{eq:rsj1}
A^T_{s}d^{(r)}+B^T_{m}d^{(j)} \leq \mu d^{(s)}, \qquad s,r,m,j=1,\ldots ,N,
\end{equation}
and inequalities \eqref{eq:theta2} and \eqref{eq:mul} hold.
Then there exists a switched diagonal L--K functional  of the form \eqref{eq:12} for system \eqref{eq:sys1}.
\end{theorem}
\textbf{Proof}  We write $P^{(s)}=\textrm{diag}\left(d^{(s)}_1/\theta_1,\ldots,d^{(s)}_n/\theta_n\right)$,
$Q_m=\tilde Q+(l-m+1)\varepsilon I$, $m=1,\ldots,l$, $s=1,\ldots,N$,
for the matrices appearing in \eqref{eq:12}.
Here $d^{(s)}_i$ and $\theta_i$ are components of the vectors $d^{(s)}$ and $\theta$ respectively, $\tilde Q$ is a positive definite diagonal matrix, and $\varepsilon>0$ is a parameter; as in the proof of Theorem \ref{thm:Swit2}, we shall show how to determine the values of $\tilde Q$, $\varepsilon$.

A relatively straightforward computation shows that the difference  of the functional \eqref{eq:12} along trajectories of \eqref{eq:sys1} satisfies the estimate
$$
\Delta \widetilde V \leq 
f^T(x(k))\left(A_{\sigma(k)}^T P^{(\sigma(k+1))} A_{\sigma(k)}    \right.
$$
$$
\left. -P^{(\sigma(k))}+  \tilde Q \right)f(x(k))
$$
$$
+2f^T(x(k)) A_{\sigma(k)}^T P^{(\sigma(k+1))} B_{\sigma(k)}f(x(k-l))
$$
$$
+
f^T(x(k-l))\left(B_{\sigma(k)}^T P^{(\sigma(k+1))} B_{\sigma(k)}-  \tilde Q \right)f(x(k-l))
$$
$$
+\varepsilon\left(l\|f(x(k))\|^2-
\sum_{m=1}^l\|f(x(k-m))\|^2\right).
$$

Consider the matrices
$$
C_{sr}= 
\left(\begin{array}{c c} 
				 A_s^T P^{(r)} A_s-P^{(s)}+  \tilde Q   &  A_s^T P^{(r)} B_s\\
			     B_s^T P^{(r)} A_s  &  B_s^T P^{(r)} B_s-  \tilde Q  			\end{array}\right)
			      			$$
for $s,r=1,\ldots,N$. We obtain
$$
C_{sr} \left(\begin{array}{c} 
				 \theta \\
			     \theta   			\end{array}\right)
\leq \left(\begin{array}{c} 
				 \lambda A_s^T d^{(r)} -d^{(s)}+\tilde Q\theta \\
			     \lambda B_s^Td^{(r)}-\tilde Q\theta   			\end{array}\right).
$$

Choose $\tilde Q$ to be the diagonal positive definite matrix such that 
$$
\tilde Q \theta=\lambda \max_{s, m=1,\ldots,N}\left\{B_s^Td^{(m)}\right\}+\delta e,  
$$
where the maximum is taken componentwise, and $\delta$ is a positive parameter.
If the value of $\delta$ is sufficiently small, then
$$
C_{sr} \left(\begin{array}{c} 
				 \theta \\
			     \theta   			\end{array}\right)\leq 
			     \left(\begin{array}{c} 
				 (\lambda\mu-1)d^{(s)}+\delta e \\
			     -\delta e   			\end{array}\right)\ll 0 
$$
for $ s,r=1,\ldots,N$.

The subsequent proof is  similar to that of Theorem \ref{thm:Swit2}. 

In the next result, we consider the case where $l=1$.  Using a virtually identical argument, we establish conditions for the existence of a closely related class of diagonal L--K functional.
\begin{proposition}\label{thm:Cor}
Consider system \eqref{eq:sys1} with $l=1$. Assume that there exist numbers $\mu>0$, $\lambda>0$ and 
vectors $d^{(1)} \gg 0,\ldots,d^{(N)}\gg 0, \theta \gg 0$ such that
\begin{equation}\label{eq:14}
A^T_{s}d^{(r)}+B^T_{r}d^{(j)} \leq \mu d^{(s)}, \qquad s,r,j=1,\ldots ,N,
\end{equation}
and that inequalities \eqref{eq:theta2} and \eqref{eq:mul} hold.
Then there exists a switched diagonal L--K functional  of the form 
\begin{eqnarray}\nonumber
\widetilde V(k,x^{(k)})=x^T(k)P^{(\sigma(k))}x(k)\\
\label{eq:15}
+f^T(x(k-1))Q^{(\sigma(k))}f(x(k-1))
\end{eqnarray}
for system \eqref{eq:sys1}, where $P^{(s)}$ and $Q^{(s)}$ are positive definite diagonal matrices, $s=1,\ldots,N$. 
\end{proposition}
\textbf{Proof} As in the proof of Theorem \ref{thm:Swit4}, we compute the difference $ \Delta \widetilde V$ along the trajectories of the system.  In this case
the difference can be bounded by quadratic forms in $f(x(k))$, $f(x(k-1))$ described by the matrices:
$$
C_{sr}= 
\left(\begin{array}{c c} 
				 A_s^T P^{(r)} A_s-P^{(s)}+   Q^{(r)}   &  A_s^T P^{(r)} B_s\\
			     B_s^T P^{(r)} A_s  &  B_s^T P^{(r)} B_s-   Q^{(s)}  			\end{array}\right)
			      			$$
for $s,r=1,\ldots,N$.  Here, the matrices $P^{(s)}$ are defined in the same way as in Theorem \ref{thm:Swit4} while the positive definite diagonal matrices $Q^{(s)}$ are given by 
$$
 Q^{(s)} \theta=\lambda \max_{m=1,\ldots,N}\left\{B_s^Td^{(m)}\right\}+\delta e,  
$$
for $s = 1, \ldots , N$. 

In a similar way as in the proof of Theorem \ref{thm:Swit4}, we can see that if $\delta>0$ is sufficiently small, then
$$
C_{sr} \left(\begin{array}{c} 
				 \theta \\
			     \theta   			\end{array}\right)\leq 
	 \left(\begin{array}{c} 
				 \lambda A_s^T d^{(r)} -d^{(s)}+ Q^{(r)}\theta \\
			     \lambda B_s^Td^{(r)}- Q^{(s)}\theta   			\end{array}\right)
			     $$
			     $$
			     \leq
			     \left(\begin{array}{c} 
				 (\lambda\mu-1)d^{(s)}+\delta e \\
			     -\delta e   			\end{array}\right)\ll 0 
$$
for $ s,r=1,\ldots,N$.

 \begin{remark}
 Theorem \ref{thm:Swit4} can be extended to systems with multiple delays.
\end{remark}

\section{Some applications of the proposed approaches}
In this section, we briefly describe some potential applications of the results and approaches outlined in the previous sections.  
\subsection{Models of digital filters}

Let the family of subsystems 
\begin{equation}\label{eq:16}
x(k+1)=f\left(A_{s} x(k)+B_{s} x(k-l)\right), \ \ s=1,\ldots,N,
\end{equation}
be given. As for the system \eqref{eq:sys1}, we assume that $x(k)\in \mathbb{R}^n$;
$A_1, \ldots, A_N$, $B_{1}, \ldots, B_{N}$ are constant nonnegative matrices; the nonlinearity  
$f: \mathbb{R}^n \rightarrow \mathbb{R}^n$ is continuous and 
diagonal, meaning:
$
f(x) = (f_1(x_1),\ldots, f_n(x_n))^T,
$
where each $f_i$ satisfies  conditions \eqref{eq:f1} and \eqref{eq:f2}; $l$ is a positive integer delay.

Systems of the form \eqref{eq:16} are used  as mathematical models of digital filters, see \cite{Kaz,EM}.

The switched
system associated with \eqref{eq:16} is 
\begin{equation}\label{eq:17}
x(k+1)=f\left(A_{\sigma(k)} x(k)+B_{\sigma(k)} x(k-l)\right).   
\end{equation}

We will first describe how to construct a common diagonal L--K functional for family \eqref{eq:16}. 
In this case, instead of \eqref{eq:V1}, we will choose such a functional in the form
\begin{eqnarray}
\nonumber
V(x^{(k)})=x^T(k)Px(k)+x^T(k-1)Q_1x(k-1) \\
\label{eq:18}
+x^T(k-2)Q_2x(k-2)
+\ldots+x^T(k-l)
Q_lx(k-l),
\end{eqnarray}
where $P,Q_1,\ldots,Q_l$ are positive definite diagonal matrices.

\begin{theorem}\label{thm:Digital-1}
If there exist real numbers $\mu>0$, $\lambda>0$ and 
vectors $d \gg 0$, $\theta \gg 0$ such that the inequalities 
\eqref{eq:d2}--\eqref{eq:mul} are valid, then there exists a common L--K functional  of the form \eqref{eq:18} for the family \eqref{eq:16}.
\end{theorem}
\textbf{Proof}  Consider the difference $\Delta V = V(x^{(k+1)}) - V(x^{(k)})$ of the functional \eqref{eq:18} along trajectories of the $s$-th subsystem from the family \eqref{eq:16} for some $s$ in $\{1,\ldots,N\}$. We obtain
$$
\Delta V =f^T\left(A_{s} x(k)+B_{s} x(k-l)\right)Pf\left(A_{s} x(k)+B_{s} x(k-l)\right)
$$
$$
+
x^T(k)(Q_1-P)x(k)
+x^T(k-1)(Q_2-Q_1)x(k-1)
$$
$$
+\ldots+x^T(k-l+1)(Q_l-Q_{l-1})x(k-l+1)
$$
$$
-x^T(k-l)
Q_lx(k-l).
$$
From \eqref{eq:f2}, it follows that
$$
\Delta V \leq \left(A_{s} x(k)+B_{s} x(k-l)\right)^T P \left(A_{s} x(k)+B_{s} x(k-l)\right)
$$
$$
+
x^T(k)(Q_1-P)x(k)
+x^T(k-1)(Q_2-Q_1)x(k-1)
$$
$$
+\ldots+x^T(k-l+1)(Q_l-Q_{l-1})x(k-l+1)
$$
$$
-x^T(k-l)
Q_lx(k-l).
$$
Taking into account this estimate, the rest of the proof follows from arguments identical to those used in the proof of Theorem \ref{thm:Swit2}.

Next, we shall present conditions of the existence of a switched diagonal L--K functional of the form 
\begin{eqnarray}
\nonumber \widetilde V(k,x^{(k)})=x^T(k)P^{(\sigma(k))}x(k)
+x^T(k-1)Q_1x(k-1)
\\ 
+x^T(k-2)Q_2 x(k-2)
\label{eq:19}
+\ldots+x^T(k-l)
Q_l x(k-l) 
\end{eqnarray}
for the system~\eqref{eq:17}. Here $P^{(1)},\ldots,P^{(N)},Q_1,\ldots,Q_l$ are positive definite diagonal matrices.

It is not too difficult to adapt the arguments of Theorems 7 and 8 to obtain the following result. 

\begin{theorem}\label{thm:Digital-2}
Assume that there exist numbers $\mu>0$, $\lambda>0$ and 
vectors $d^{(1)} \gg 0,\ldots,d^{(N)}\gg 0, \theta \gg 0$ such that the inequalities
\eqref{eq:theta2}, \eqref{eq:mul} and \eqref{eq:rsj1} are satisfied.
Then there exists a switched diagonal L--K functional  of the form \eqref{eq:19} for system \eqref{eq:17}.
\end{theorem}

\subsection{Discrete-time neural networks}

Next, consider the switched system
\begin{equation}\label{eq:20}
x(k+1)=A_{\sigma(k)} f(x(k))+B_{\sigma(k)} f(x(k-l))+u(k)   
\end{equation}
and the associated family of subsystems
\begin{equation}\label{eq:21}
x(k+1)=A_{s} f(x(k))+B_{s} f(x(k-l))+u(k),
\end{equation}
$$
s=1,\ldots,N.
$$
Here $u(k)$ is a bounded input, and the rest of the notation is the same as for the original system given in \eqref{eq:sys1}.
Systems of the form \eqref{eq:20} arise in models of neural networks, see \cite{Kaz,TG} for background on their use in this context.

In general, the family \eqref{eq:21} will not admit a common equilibrium position. However, the approaches outlined in this paper permit us to derive conditions under which solutions of \eqref{eq:20} will be ultimately bounded. 

Let $x\left(k,x^{(k_0)},k_0\right)$ denote the solution of \eqref{eq:20}
with initial conditions $k_0\geq 0$, $x^{(k_0)}\in \mathbb{R}^{(l+1)n}$.

\begin{definition}
The system \eqref{eq:20} is said to be uniformly ultimately bounded
with the ultimate bound $R$ 
if, for any $D>0$, there exists a positive integer $\tilde k$ such 
that $\left\|x\left(k,x^{(k_0)},k_0\right)\right\| \leq  R$ for all $k_0\geq 0$, $\left\|x^{(k_0)}\right\| \leq D$, $k\geq \tilde k+  k_0$ and for arbitrary switching law $\sigma(k)$.
\end{definition}

\begin{theorem}\label{thm:Networks-1}
Let $f_i(x_i)\to -\infty$ as $x_i\to -\infty$, and
$f_i(x_i)\to +\infty$ as $x_i\to +\infty$, $i=1,\ldots,n$. Assume that  there exist numbers $\mu>0$, $\lambda>0$ and 
vectors $d \gg 0$, $\theta \gg 0$ such that inequalities 
\eqref{eq:d2}--\eqref{eq:mul} are valid. Then there exists a common L--K functional  of the form \eqref{eq:V1} for the family \eqref{eq:21} guaranteeing uniform ultimate boundedness of system \eqref{eq:20}.
\end{theorem}
\textbf{Proof} Define matrices $P,Q_1,\ldots,Q_l$
in the same way as in the proof of Theorem~4. Consider the difference of  the functional \eqref{eq:V1} along trajectories of the $s$-th subsystem from the family \eqref{eq:21} for some $s$ in $\{1,\ldots,N\}$. For sufficiently small values of the parameters $\varepsilon$ and $\delta$, we obtain 
\begin{equation}\label{eq:22}
\Delta V \leq -\beta_1\sum_{j=0}^l  \|f(x(k-j))\|^2+\beta_2,
\end{equation}
where $\beta_1$ and $\beta_2$ are positive constants.

Using the estimate \eqref{eq:22}, the remainder of the proof is very similar to the proof of Theorem 5.2 in \cite{Al-Diff}.

\begin{theorem}\label{thm:Networks-2}
Let $f_i(x_i)\to -\infty$ as $x_i\to -\infty$, and
$f_i(x_i)\to +\infty$ as $x_i\to +\infty$, $i=1,\ldots,n$.
Assume that there exist numbers $\mu>0$, $\lambda>0$ and 
vectors $d^{(1)} \gg 0,\ldots,d^{(N)}\gg 0, \theta \gg 0$ such that inequalities
\eqref{eq:theta2}, \eqref{eq:mul} and \eqref{eq:rsj1} are valid.
Then there exists a switched diagonal L--K functional  of the form \eqref{eq:12} guaranteeing uniform ultimate boundedness of system \eqref{eq:20}.
\end{theorem}

The proof of the theorem is similar to that of Theorem \ref{thm:Swit4}.

\section{Numerical example}\label{sec4}
In this section, we present a simple numerical example in order to illustrate how our results relate to each other and to earlier work.

Consider a family \eqref{eq:sys2} consisting of two subsystems of dimension 2 ($N=2$, $n=2$) with the system matrices
$$
A_1 =\frac14 \left(\begin{array}{c c} 
				 0  &  0\\
			     1  &  1  
			\end{array}\right),  \qquad
B_1 = \frac14\left(\begin{array}{c c}
		a &  1\\
	    2 &  0
	    \end{array}\right),
$$
$$
A_2 =\frac14 \left(\begin{array}{c c} 
				 0  &  1\\
			     0  &  2  
			\end{array}\right),  \qquad
B_2 = \frac14\left(\begin{array}{c c}
		1 &  1\\
	    0 &  0
	    \end{array}\right),
$$
             where $a$ is a positive parameter.
In addition, assume that $l=1$.

In this case, the system of inequalities \eqref{eq:d1} admits a positive solution if and only if $a<1$, whereas \eqref{eq:theta1} admits a positive solution if and only if $a<3$. Thus, we can apply Theorem 2 only for $a<1$.

Next, let $a=2$. For this value of the parameter $a$,  \eqref{eq:theta2} admits a positive solution if and only if 
$\lambda\geq (1+\sqrt{6})/4$, whereas
system \eqref{eq:d2} admits a positive solution if and only if $\mu\geq (2+\sqrt{6})/4$.
Verifying condition \eqref{eq:mul}, we obtain 
$$
\frac{2+\sqrt{6}}4 \cdot  \frac{1+\sqrt{6}}4 <1.
$$
Thus, Theorem 4 guarantees the existence of a common L--K functional of the form \eqref{eq:V1}.

Finally,  choose $a=9/4$. Then system \eqref{eq:theta2} admits a positive solution if and only if 
$\lambda\geq (13+\sqrt{217})/32$, and
system \eqref{eq:d2} admits a positive solution if and only if $\mu\geq (17+\sqrt{433})/32$.

Since 
$$
\frac{13+\sqrt{217}}{32} \cdot  \frac{17+\sqrt{433}}{32} >1,
$$
 condition \eqref{eq:mul} is not fulfilled. Hence, for this value of parameter $a$,we can not apply Theorem~4.

On the other hand, if we consider system \eqref{eq:14} for $a=9/4$, $\mu=1.152$, then the system admits the solution 
$d^{(1)}=(1.179, 0.5)^T$, $d^{(2)}=(1.3, 1)^T$. 
We obtain
$$
\frac{13+\sqrt{217}}{32} \cdot  1.152<1.
$$
Hence, all the conditions  of Proposition~4 are fulfilled, and we can  guarantee the existence of a switched diagonal L--K functional of the form \eqref{eq:15} for the associated switched system.

\section{Conclusion}\label{sec5}
In this short paper, we have described sufficient conditions for the existence of common and switched diagonal Lyapunov functionals for switched positive nonlinear systems subject to time-delay.  In particular, the result of Theorem \ref{thm:Swit2} relaxes the condition required in \cite{AlexMas} for the existence of a common diagonal Lyapunov functional, thereby giving a less conservative stability criterion.  It also provides an extension of Theorem 3 of \cite{Past} to nonlinear time-delayed systems.  The result of Theorem \ref{thm:Swit4} describes conditions for a switched diagonal Lyapunov functional to exist for the same class of nonlinear switched systems, and gives an alternative type of condition to the LMIs described in \cite{Shaker}.

\section{Acknowledgments}

This work was partially supported by the Government of the Russian Federation (Grant no.  074-U01),
the Russian Foundation for Basic Research (Grant no. 16-01-00587), 
the Science Foundation Ireland (Grant no. 13/RC/2094), and the European Regional Development Fund through the Southern \& Eastern Regional Operational Programme to Lero -- the Irish Software Research Centre (www.lero.ie).


\begin{thebibliography}{99}







\bibitem{ValCol}
Blanchini F., Colaneri P., Valcher M.E.: 
`Switched positive linear systems', 
        \textit{Foundations and Trends in Systems and Control}, 2015, \textbf{2}, no.~2, 
        pp.~101--273
        
\bibitem{FarRin}
 Farina L.,  Rinaldi S.:
`Positive linear systems: theory and applications'
(Wiley, New York, 2000)

\bibitem{ShoTCP}
Shorten R.N., Wirth F., Leith D.: `A positive systems model of TCP-like congestion control', \textit{IEEE Trans. Networking},
2006, \textbf{14}, no.~3, pp.~616--629

\bibitem{ZhaoLiuRees09}
Zhao Y.B., Liu G.P., Rees D.: `Stability and stabilisation of discrete-time
networked control systems: a new time
delay system approach', \textit{IET Control Theory Appl.},
2009, \textbf{4}, no.~9, pp.~1859--1866



\bibitem{MasSho07}
Mason O., Shorten R.N.: `On linear copositive Lyapunov functions and the stability of switched positive linear systems', \textit{IEEE Trans. Aut. Cont.},
2007, \textbf{52}, no.~7, pp.~1346--1349

\bibitem{Kaz} 
 Kaszkurewicz E.,  Bhaya A.: `Matrix diagonal stability in systems and computation' (Birkh\"{a}user, Boston, Basel, Berlin, 2000)
 

\bibitem{Lam}
 Wu L.,  Lam J.,  Shu Z.,~et al.: `On stability and stabilizability of positive delay systems', \textit{Asian Journal of Control}, 2009, \textbf{11}, pp.~226--234.
 



\bibitem{AlexMas} 
Aleksandrov A.,  Mason O.:
`Diagonal Lyapunov--Krasovskii functionals for dis\-cre\-te-time positive systems with delay',
\textit{Syst. Control Lett.}, 2014, \textbf{63}, pp.~63--67



\bibitem{HJ}
 Horn R.A.,  Johnson C.R.:
`Matrix analysis'
(Cambridge University Press, New York, 1985)



\bibitem{Short}
Shorten R.,  Wirth F., Mason O.,~et al.: `Stability criteria for switched and hybrid systems', \textit{SIAM Rev.}, 2007, \textbf{49},  no. 4,
pp.~545--592



\bibitem{Past}
Pastravanu O.C., Matcovschi M.-H.: 
`Max-type copositive Lyapunov functions for switching positive linear systems', 
        \textit{Automatica}, 2014, \textbf{50}, 
        pp.~3323--3327



\bibitem{Huang} 
 Zhang J.,  Huang J.,  Zhao X.:
  `Further results on stability and stabilisation of switched positive systems', \textit{IET Control Theory Appl.},
2015, \textbf{9}, no.~14, pp.~2132--2139
  
  

\bibitem{AiT} 
Aleksandrov A.,  Platonov A.:
`On absolute stability of
one class of nonlinear switched systems',
\textit{Autom. Remote
Control}, 2008, \textbf{69}, no.~7, pp.~1101--1116



\bibitem{Shaker}
 Shaker H.R., How J.P.: 
'Stability analysis for class of switched nonlinear systems'. Proc.
    American Control Conf.,
Marriott Waterfront, Baltimore, MD, USA,
June 30-July 02 2010,  pp.~2517--2520




\bibitem{ZhenFen12}
Zheng Y., Feng G.: `Diagonal stabilisation of a class of single-input
discrete-time switched systems', \textit{IET Control Theory Appl.},
2012, \textbf{7}, no.~4, pp.~515--522





\bibitem{Knorn}
Knorn F., Mason O., Shorten R.N.: 
`On linear co-positive Lyapunov functions for sets of linear positive systems', 
        \textit{Automatica}, 2009, \textbf{45}, 
        pp.~1943--1947




\bibitem{ForVal10}
Fornasini E., Valcher M.E..: `Linear copositive functions for continuous time positive switched systems', \textit{IEEE Trans. Aut. Cont.},
2010, \textbf{55}, no.~8, pp.~1933--1937





\bibitem{Doan}
Doan T.S., Kalauch A., Siegmund S.: `A constructive approach to linear Lyapunov functions for positive switched systems using Collatz-Wielandt sets', \textit{IEEE Trans. Aut. Cont.},
2013, \textbf{58}, no.~3, pp.~748--751

\bibitem{Gowda} 
Song Y., Seetharama-Gowda M., Ravindran G.:
`On some properties of P-matrix sets',
\textit{Lin. Alg. and Appl.}, 1999, \textbf{290}, pp.~237--246


 




\bibitem{Phi}
Sun X., Li J., Zhao J.: `Stabilization for a class of discrete-time switched $\Phi$-systems', \textit{Circuits, Systems, and Signal Processing}, 2017, \textbf{36}, pp.~834--844

\bibitem{Sontag}  Sontag E.D., Karny M., Warwick K.,~et al.: `Recurrent neural networks: some systems-theoretic aspects in dealing with complexity: a neural network approach'
(Springer, London, 1997)


\bibitem{Liao}
 Liao X.,  Yu P.:
`Absolute stability of nonlinear control systems'
(Springer, New York, Heidelberg, 2008)


\bibitem{EM}
Erickson K.T., Michel A.N.: `Stability analysis of fixed-point digital filters using computer generated  Lyapunov functions -- Part I: Direct form and coupled form filters', \textit{IEEE Trans. on Circuits and Systems},
1985, \textbf{32}, pp.~113--132


\bibitem{Kam}
 Kamenetskiy V.A.,  Pyatnitskiy Ye.S.: `An iterative method of Lyapunov function construction for differential inclusions', 
 \textit{Syst. Control Lett.}, 1987, \textbf{8}, pp.~445--451


\bibitem{Boid}  Boyd S., Ghaoui E., Feron E.,~et al.: `Linear matrix inequalities in system and control theory'
(SIAM, Philadelphia, 1994)






























\bibitem{Liu}
Liu X.: `Stability analysis of switched positive systems: a switched linear copositive Lyapunov function method', \textit{IEEE Trans. Circuits and Systems} (II), 2009,
\textbf{56}, pp.~414--418


 


\bibitem{Daaf}
Daafous J., Riedinger P., Iung C.: `Stability analysis and control synthesis for switched systems: a switched Lyapunov function approach', \textit{IEEE Trans. Autom. Control},
2002, \textbf{47}, pp.~1883--1887
















\bibitem{zapp}
 Zappavigna A.,  Colaneri P.,  Geromel J.C.,~et al.:
'Dwell time analysis for continuous-time switched linear positive systems'. Proc.
    American Control Conf.,
Marriott Waterfront, Baltimore, MD, USA,
June 30-July 02 2010,  pp.~6256--6261










\bibitem{TG}
Takeda M., Goodman J.W.: `Neural networks for computation: number representations and programming complexity', \textit{Applied Optics},
1986, \textbf{25}, pp.~3033--3046










\bibitem{Al-Diff} 
Aleksandrov A.,  Chen Y., Platonov A.,~et al.:
`Stability analysis and uniform
ultimate boundedness control synthesis for a class of nonlinear switched difference systems',
\textit{J. Difference Equ. Appl.}, 2012, \textbf{18}, pp.~1545--1561



\end{thebibliography}
\end{document}